\documentclass[11pt]{amsart}

\usepackage[T1]{fontenc}
\usepackage{lmodern}
\usepackage{microtype}
\usepackage{amsmath,amssymb,mathtools}
\usepackage{enumitem}
\usepackage{xcolor}
\usepackage[margin=1.05in]{geometry}
\usepackage[colorlinks=true,linkcolor=blue!55!black,citecolor=blue!55!black,urlcolor=blue!55!black]{hyperref}
\usepackage{tikz}
\usepackage{ulem}

\newtheorem{theorem}{Theorem}[section]
\newtheorem{proposition}[theorem]{Proposition}
\newtheorem{corollary}[theorem]{Corollary}
\newtheorem{lemma}[theorem]{Lemma}
\theoremstyle{definition}

\newtheorem{example}[theorem]{Example}
\theoremstyle{remark}
\newtheorem{remark}[theorem]{Remark}

\newcommand{\R}{\mathbb{R}}

\newcommand{\eps}{\varepsilon}
\newcommand{\vol}{\operatorname{area}}

\newcommand{\cL}{\mathcal{L}}
\newcommand{\cF}{f}

\def\srem#1{\ifvmode\leavevmode\fi\vadjust{\vbox to0pt{\vss 
    \hbox to 0pt{\hskip\hsize\hskip1em           
 \vbox{\hsize2.2cm\tiny\raggedright\pretolerance10000
 \noindent #1\hfill}\hss}\vbox to8pt{\vfil}\vss}}}%

\title[Floating Body Asymptotics via Fractal Strings]
{Floating Body Asymptotics via Fractal Strings}

\author{Christoph Th\"ale}
\address{Faculty of Mathematics,
Ruhr University Bochum,
44780 Bochum, Germany}
\email{christoph.thaele@rub.de}

\author{Steffen Winter}
\address{Institute of Stochastics,
Karlsruhe Institute of Technology,
76131 Karlsruhe, Germany}
\email{steffen.winter@kit.edu}

\date{}


\begin{document}

\begin{abstract}
Sch\"utt and Werner associated an almost polygonal planar convex body with a decreasing convex sequence and obtained two-sided estimates for the area removed by its floating body. We reinterpret the gaps of the defining sequence as the lengths of a normalized fractal string. Our main result shows that the area deficit, that is, the difference between the area of the set and its floating body, is comparable to an integral transform of the inner tubular volume of this string. As a consequence, if the  string has Minkowski dimension \(D\), then the principal power exponent of the area deficit is \(1-D/3\). We also derive refinements for regularly varying tube functions, establish a realization criterion extending the almost polygonal case of the theorem of Sch\"utt and Werner, and discuss several examples. As an application, we obtain corresponding asymptotic rates for the expectations and variances of the missed area and the number of vertices of random polygons, including new intermediate behaviors between the classical polygonal and smooth cases.
\end{abstract}

\maketitle

\section{Introduction}

Let \(K\subset\R^2\)
be a convex body. For \(\delta>0\), its convex floating body \(K_\delta\)
is obtained by intersecting all closed half-planes whose complementary
cap in \(K\) has area at most \(\delta\). The resulting area deficit $\Delta_K(\delta):=\vol(K\setminus K_\delta)$
measures the size of the boundary layer detected at cap scale
\(\delta\). Floating bodies arise naturally in affine convex geometry and in the
approximation of convex bodies by polytopes, see, for example,
\cite[Section~10.6]{Schneider2014} and the survey article \cite{Werner2022}. In the smooth setting, the
first-order behavior of \(\Delta_K(\delta)\) as $\delta\downarrow 0$ is governed by the classical affine
perimeter. More precisely, if \(\partial K\) is of class \(C^2\) and has
everywhere positive curvature $\kappa(x)>0$, $x\in\partial K$, then
\[
\Delta_K(\delta)
\asymp
\delta^{2/3},
\qquad \delta\downarrow0,
\]
and the normalized area deficit has a limit given by a multiple of the affine perimeter $\int_{\partial K}\kappa(x)^{1/3}\,ds(x)$ in which \(ds\) denotes the arc-length element on \(\partial K\), see \cite{SchuettWerner1990}. Here and throughout, \(f\asymp g\) means that the ratio \(f/g\)
is bounded above and below by positive constants in the asymptotic
regime under consideration. At the opposite end of the
regularity spectrum, if \(K\) is a polygon, then
\[
\Delta_K(\delta)
\asymp
\delta\log\frac1\delta,
\qquad \delta\downarrow0,
\]
see \cite{Schuett1991} and
\cite[Corollary~1.2]{BesauSchuettWerner2018} for the
higher-dimensional formulation in terms of flag numbers. Again, the normalized limit exists and is given by {\color{red}\sout{a}} $1/2$ times the number of vertices of the polygon.  These two cases exhibit fundamentally different
mechanisms. Near a smooth boundary point, a cap of area \(\delta\) has
tangential diameter of order \(\delta^{1/3}\), whereas near a vertex of
a polygon the accumulation of contributions over a range of cap scales
produces the logarithmic factor.

\begin{figure}[t]
	\centering
	\begin{minipage}{0.48\textwidth}
		\centering
		\includegraphics[
		trim={0 0 0 10mm},
		clip,
		width=\textwidth
		]{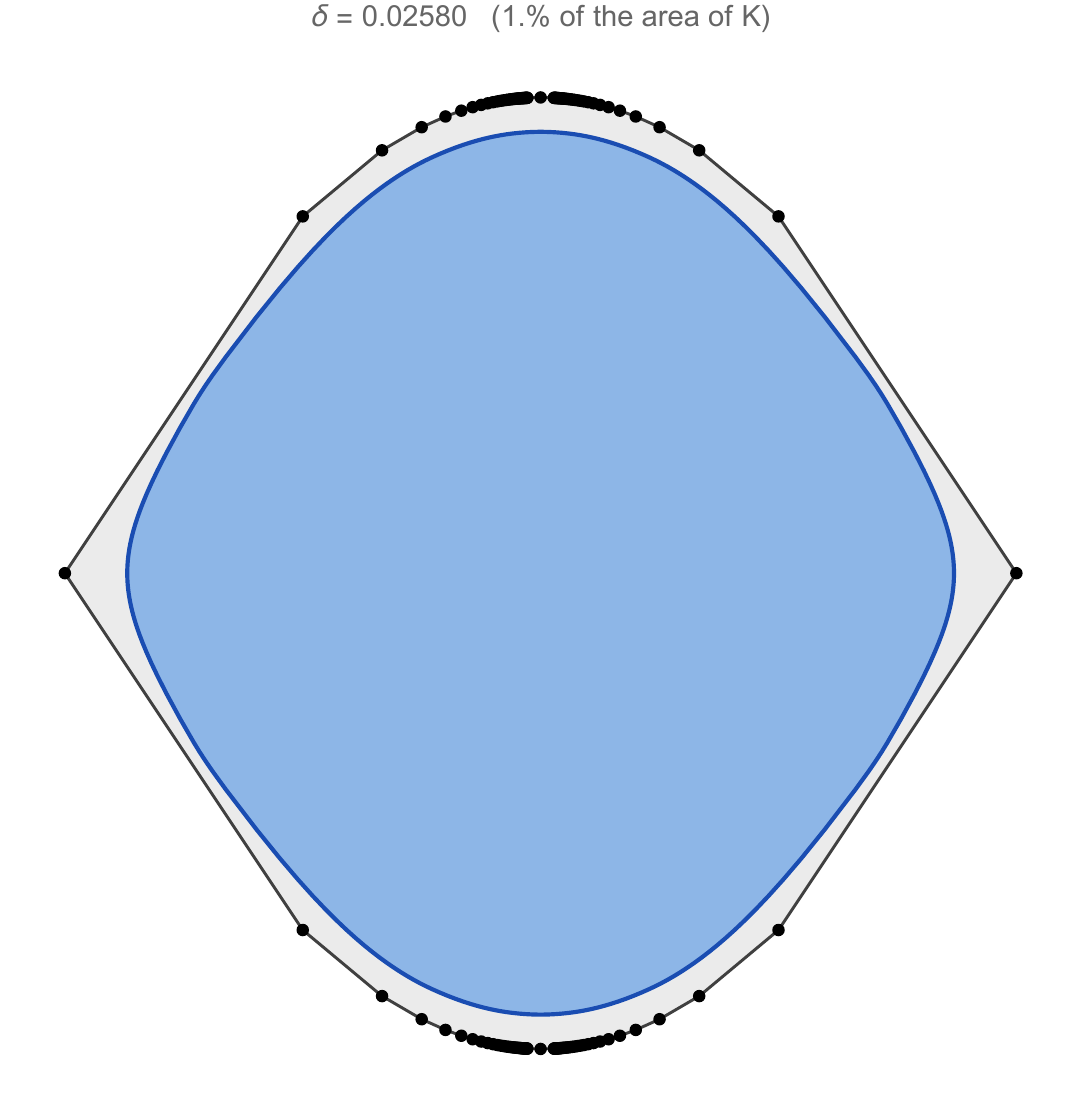}
		
		\smallskip
		\(a_n=\frac{1}{n+1}\)
	\end{minipage}
	\hfill
	\begin{minipage}{0.48\textwidth}
		\centering
		\includegraphics[
		trim={0 0 0 10mm},
		clip,
		width=\textwidth
		]{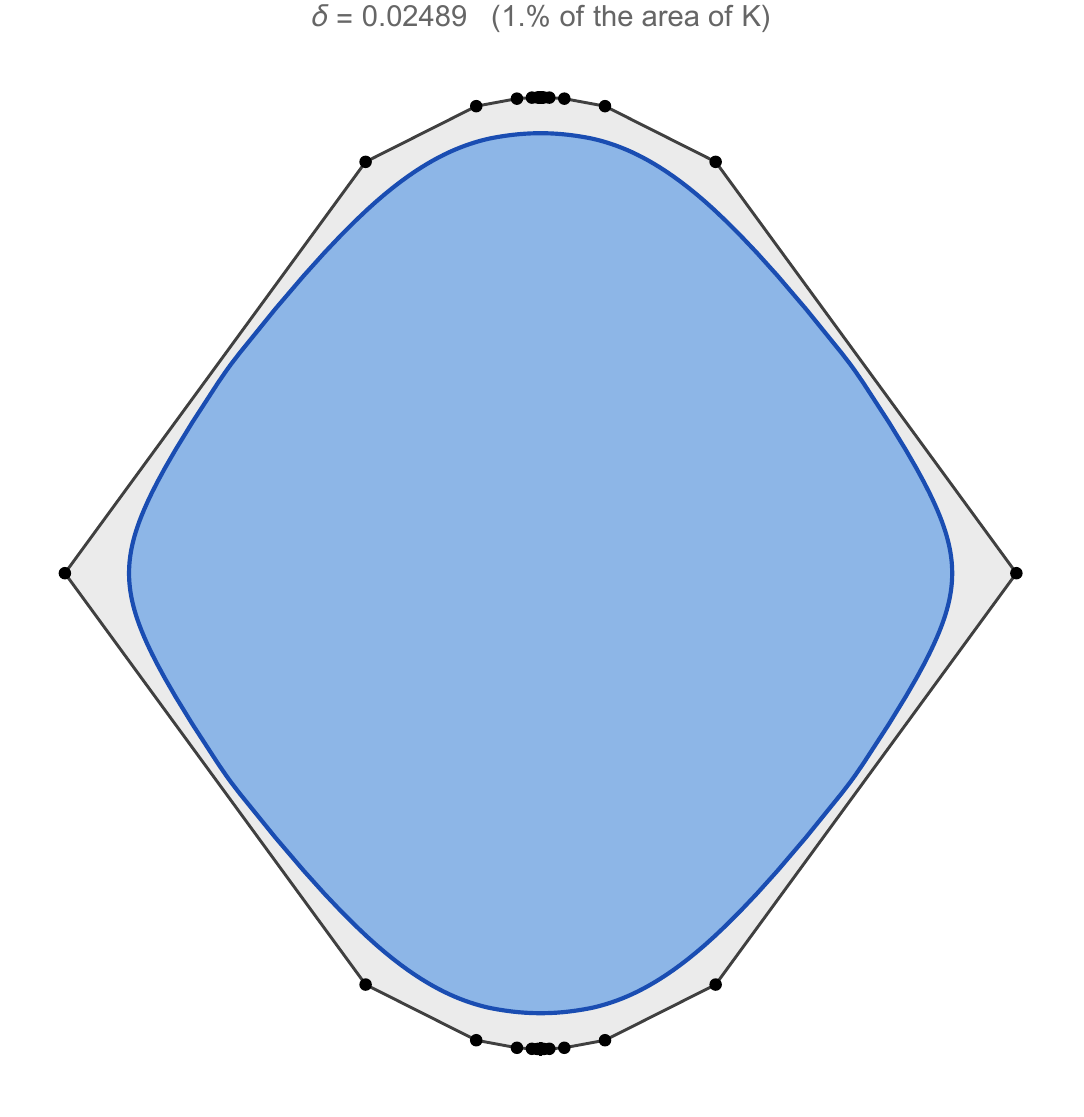}
		
		\smallskip
		\(a_n=e^{-n}\)
	\end{minipage}
	
	\caption{Almost polygonal bodies (grey) and their floating
		bodies (blue) for a polynomially and an exponentially
		decreasing tail sequence.}
	\label{fig:floating-bodies-examples}
\end{figure}

For a general nonsmooth and nonpolygonal convex body, the small-cap
geometry may vary substantially along the boundary, and neither of the
two preceding asymptotic regimes need describe the true order of the
area deficit. A tractable intermediate class was introduced by
Sch\"utt and Werner in \cite{SchuettWerner1992}. Starting from a decreasing convex sequence $1=a_0>a_1>a_2>\ldots$ converging to zero, they constructed an almost polygonal convex body whose extreme points
lie on two parabolas and accumulate at two boundary points, see Figure \ref{fig:floating-bodies-examples}. For this
class they obtained explicit two-sided estimates for
\(\Delta_K(\delta)\) and proved a realization theorem showing
that a wide range of intermediate rates can occur.

The purpose of this article is to show that the gap sequence in this
construction is naturally encoded by a fractal string and that the tube geometry of this string determines the floating-body asymptotics. Put
\(\ell_n:=a_n-a_{n+1}\), \(n\geq0\). Since the sequence \((a_n)_{n\geq 0}\) is decreasing and
convex, and converges to $0$,
\[
\ell_0\geq\ell_1\geq\cdots>0,
\qquad
\sum_{n=0}^{\infty}\ell_n=1.
\]
Thus \(\cL=(\ell_n)_{n\geq0}\) is a normalized fractal string in the
sense of \cite{LapidusVanFrankenhuijsen2000}. For this string we use the canonical realization as a subset of $\R$ obtained by placing
intervals of lengths \(\ell_0,\ell_1,\ldots\) consecutively in
\((0,1)\). Equivalently, if
\(a_n=\sum_{j=n}^{\infty}\ell_j\), then its boundary is
\[
A_{\cL}:=\{0\}\cup\{a_n:n\geq0\}.
\]
We denote by \(V_{\cL}(r)\) the corresponding inner tubular volume. All these notions will be recalled in Section~\ref{sec:Prelim}.

Our main result identifies the geometric quantity that controls the
area deficit. Namely, if \(K(\cL)\) denotes the almost polygonal body
associated with a normalized fractal string \(\cL\), we show that
\begin{equation*}
	\Delta_{K(\cL)}(\eps^3)
	\asymp
	\eps^3
	\int_\eps^1
	\frac{V_{\cL}(r)}{r^2}\,dr,
	\qquad \eps\downarrow0.
\end{equation*}
This representation gives a direct relation between the fractal
geometry of the string and the floating-body asymptotics. Namely, if the
Minkowski dimension $D=\dim_M A_{\cL}$  of the string exists, then
\begin{equation}\label{eq:intro-dimension}
	\Delta_{K(\cL)}(\delta)
	=
	\delta^{1-D/3+o(1)},
	\qquad \delta\downarrow0.
\end{equation}
The dimension alone does not, however, detect logarithmic or more
general slowly varying corrections. But if $V_{\cL}(r)\asymp r^{1-D}L(1/r)$,
where \(L\) is slowly varying and \(D>0\), then we show that
\[
\Delta_{K(\cL)}(\delta)
\asymp
\delta^{1-D/3}L(\delta^{-1/3}),
\]
whereas at the critical value \(D=0\), an additional logarithmic integration
appears.  We also establish a realization theorem for prescribed floating-body
deficit rates by inverting the integral transform of the tube
function. In the almost polygonal case, this extends the realization
theorem of \cite{SchuettWerner1992} under weaker regularity and
growth assumptions.

There are also direct probabilistic consequences of our results. To describe them, let
\(X_1,X_2,\ldots\) be independent random points uniformly distributed
in a planar convex body \(K\), and let $P_n:=\operatorname{conv}(X_1,\ldots,X_n)$, $n\geq 3$,
be the random polygon in $K$ generated by the first $n$ points. A theorem of
B\'ar\'any and Larman \cite{BaranyLarman1988} compares the expected missed area with the area deficit. It says that
\begin{equation*}
	\mathbb E\vol(K\setminus P_n)
	\asymp
	\Delta_K\left(\frac{\vol(K)}{n}\right),
	\qquad n\to\infty.
\end{equation*}
Consequently, for random polygons $P_n$ in an almost polygonal body $K(\cL)$ associated with a normalized fractal string $\cL$, relation \eqref{eq:intro-dimension} implies that
\[
\mathbb E\vol\bigl(K(\cL)\setminus P_n\bigr)
=
n^{-1+D/3+o(1)}.
\]
Moreover, the asymptotics for the expected number of vertices of $P_n$ follow from Efron's identity. Writing
\(f_0(P_n)\) for the number of vertices of the random polygons $P_n$ in almost polygonal bodies considered here, we obtain
\[
\mathbb E f_0(P_n)
=
n^{D/3+o(1)}.
\]
In other words, the Minkowski dimension of the fractal string $\cL$
controls both the decay of the expected missed area and
the growth of the expected combinatorial complexity of the associated
random polygon. Varying the fractal string will allow us to produce asymptotic regimes that
interpolate between the classical polygonal and smooth cases, including
nontrivial iterated logarithmic corrections. Combining our results with the variance estimates of Pardon
\cite{Pardon2012}, we also determine the variance orders and hence
the fluctuation scales in the corresponding central limit theorems
for random polygons.

Both convex floating bodies and random polytopes have been studied
extensively in \(\R^d\) for arbitrary \(d\geq2\), see, for example,
\cite{BaranyLarman1988,Reitzner2010,Schneider2014,SchuettWerner1990}. Our results, however, are specific to the planar setting, where the geometry of almost polygonal bodies admits a natural description in terms of fractal strings.

\medspace

The remaining parts are organized as follows. In
Section~\ref{sec:Prelim} we recall the almost-polygonal construction
of Sch\"utt and Werner and introduce the necessary material on fractal
strings, tube volumes, and Minkowski dimension. In
Section~\ref{sec:Results} we formulate and prove our main results. Section~\ref{sec:Ex}
contains a collection of examples, and the final Section~\ref{sec:random-polygons}
discusses the consequences for random polygons.

\section{Preliminaries}\label{sec:Prelim}

\subsection{Floating bodies and the almost-polygonal body construction}

For a unit vector \(u\in S^1\) and \(t\in\R\), write $H^-(u,t)=\{x\in\R^2:\langle x,u\rangle\le t\}$,
where \(\langle\cdot,\cdot\rangle\) denotes the standard scalar product
on \(\R^2\), and let \(\vol(\cdot)\) denote Lebesgue measure on \(\R^2\).
For a convex body $K\subset\R^2$ the (convex) floating body is defined as
\[
    K_\delta
    =\bigcap\bigl\{H^-(u,t):
       \vol\bigl(K\setminus H^-(u,t)\bigr)\le\delta\bigr\},
\]
see \cite{SchuettWerner1990}. In other words, \(K_\delta\) is obtained from \(K\) by removing all boundary caps of area at most \(\delta\). The area deficit between $K$ and its $\delta$-floating body $K_\delta$ is denoted by
\[
\Delta_K(\delta):=\vol(K\setminus K_\delta).
\]

We now recall the construction of the almost polygonal bodies from
\cite{SchuettWerner1992}. Consider a sequence
\(a=(a_n)_{n\geq0}\) of non-negative real numbers satisfying
\begin{equation}\label{eq:a-assumptions}
    a_0=1,\qquad a_n\downarrow0,
    \qquad
    a_{n-1}-a_n\ge a_n-a_{n+1}>0,\qquad n\geq 1.
\end{equation}
Following \cite{SchuettWerner1992}, we define the almost polygonal body $K(a)$ as
\begin{equation}\label{eq:Ka}
	K(a)
	:=
	\operatorname{conv}
	\left(
	\bigl\{
	(\sigma a_n,\tau(a_n^2-1)):
	n\geq0,\ \sigma,\tau\in\{-1,1\}
	\bigr\}
	\cup\{(0,-1),(0,1)\}
	\right).
\end{equation}
The extreme points of \(K(a)\) lie on the two parabolas with equations
\(y=x^2-1\) and \(y=-x^2+1\), and their only accumulation points are
\((0,-1)\) and \((0,1)\), see Figure~\ref{fig:floating-bodies-examples}.

The gap lengths determined by the sequence $a=(a_n)_{n\geq 0}$ are denoted by $\ell_n:=a_n-a_{n+1}$, $n\geq 0$, and, for $\delta>0$, let $n(\delta)$ be the smallest integer for which
\begin{equation}\label{eq:ndef}
    \ell_{n(\delta)}\le 2^{1/3}\delta^{1/3}.
\end{equation}
The estimate that drives the present article is the following reformulation of Proposition~2.1 in \cite{SchuettWerner1992}. It describes the area deficit between $K(a)$ and its $\delta$-floating body in terms of the gap sequence $(\ell_n)_{n\geq 0}$ and the number $n(\delta)$.

\begin{proposition}\label{prop:SW}
Let \(a=(a_n)_{n\geq0}\) satisfy \eqref{eq:a-assumptions}, let
\(K(a)\) be defined by \eqref{eq:Ka}, put
\(\ell_n=a_n-a_{n+1}\), and let \(n(\delta)\) be defined by
\eqref{eq:ndef}. Then there are constants \(c\geq1\) and
\(\delta_0>0\) such that, for \(0<\delta<\delta_0\),
\begin{align*}
 &\frac1c\,\delta^{2/3}a_{n(\delta)}
 +\frac1c\,\delta\sum_{k=2}^{n(\delta)}
 \log\!\left(
 \frac{\ell_{k-1}\ell_{k-2}(\ell_{k-1}+\ell_{k-2})}{2\delta}
 \right)
 \le \Delta_{K(a)}(\delta)  \\
 &\hspace{1cm}\le
 c\,\delta^{2/3}a_{n(\delta)}+c\,\delta n(\delta)
 +c\,\delta\sum_{k=2}^{n(\delta)}
 \log\!\left(
 \frac{\ell_{k-1}\ell_{k-2}(\ell_{k-1}+\ell_{k-2})}{2\delta}
 \right).
\end{align*}
\end{proposition}

The two contributions in Proposition~\ref{prop:SW} reflect different
types of boundary behavior. Near the accumulation points of the
extreme points, the gaps between consecutive vertices become too
small to be resolved at cap-area scale \(\delta\). The boundary then
behaves approximately like the underlying parabolas, giving rise to
the contribution \(\delta^{2/3}a_{n(\delta)}\). The logarithmic sum
accounts for the contributions from the remaining vertices, where
the boundary exhibits polygonal behavior.
The transition between these two regions occurs when the gap lengths
are of order \(\delta^{1/3}\). Indeed, a cap of area \(\delta\) cut
from a parabola has height of order \(\delta^{2/3}\), while its
projection onto the horizontal axis has length of order
\(\delta^{1/3}\). This explains the scale used in
\eqref{eq:ndef}.

\begin{remark}
As we shall see in the proof of Theorem~\ref{thm:main},
the second term \(\delta n(\delta)\) in the upper bound of
Proposition~\ref{prop:SW} is always dominated by the
logarithmic sum and can therefore be omitted after
adjusting the constant $c$.
\end{remark}

\subsection{Fractal strings, tube volumes, and Minkowski dimension}

By a normalized fractal string we understand a nonincreasing sequence
\[
\cL=(\ell_n)_{n\geq0},
\qquad
\ell_0\geq\ell_1\geq\cdots>0,
\qquad
\sum_{n=0}^{\infty}\ell_n=1.
\]
Such a sequence may be realized as the sequence of lengths of the
connected components of a bounded open subset of the real line. For
our purposes, the following canonical realization is particularly
convenient. Define the tail sums
\begin{equation}\label{eq:tails}
	a_n:=\sum_{j=n}^{\infty}\ell_j,
	\qquad n\geq0,
\end{equation}
and put
\[
A_{\cL}:=\{0\}\cup\{a_n:n\geq0\},
\qquad
\Omega_{\cL}:=(0,1)\setminus A_{\cL}.
\]
Since \(a_0=1\), \(a_n\downarrow0\), and $a_n-a_{n+1}=\ell_n$,
we have the disjoint decomposition
\[
\Omega_{\cL}
=
\bigsqcup_{n=0}^{\infty}(a_{n+1},a_n).
\]
In particular, the connected components of \(\Omega_{\cL}\) have
lengths \(\ell_n\), and $\partial\Omega_{\cL}=A_{\cL}$. Conversely, if a sequence \(a=(a_n)_{n\geq0}\) satisfies
\eqref{eq:a-assumptions}, then $\ell_n:=a_n-a_{n+1}$, $n\geq 0$,
defines a normalized fractal string $\cL(a)$ written in nonincreasing order, and
\(a\) is recovered from \(\cL(a)\) through \eqref{eq:tails}.

The inner tubular volume of \(\cL=(\ell_n)_{n\geq0}\) is defined by
\[
V_{\cL}(r)
:=
\bigl|
\{x\in\Omega_{\cL}:
\operatorname{dist}(x,A_{\cL})<r\}
\bigr|,
\qquad r>0,
\]
where $\operatorname{dist}(x,A)$ denotes the
Euclidean distance from a point \(x\) to a set \(A\subset\R\), and \(|\cdot|\)
denotes Lebesgue measure on \(\R\). A connected component of $\Omega_{\cL}$
of length $\ell_n$ contributes precisely
\(\min\{2r,\ell_n\}\) to the inner tubular volume $V_{\cL}(r)$. Consequently,
\begin{equation}\label{eq:tube-sum}
	V_{\cL}(r)
	=
	\sum_{n=0}^{\infty}\min\{2r,\ell_n\},
\end{equation}
see also \cite[Equation~(1.9)]{LapidusVanFrankenhuijsen2000}.

Let \(A\subset[0,1]\) be a compact set, and let \(N_A(r)\) denote the smallest
number of intervals of length \(r>0\) required to cover \(A\). The following comparison between covering numbers and
parallel volumes will be useful, see also the proof of \cite[Proposition~3.2]{Falconer2003}.

\begin{lemma}\label{lem:tube-cover}
	For every compact set \(A\subset[0,1]\) and every \(r\in(0,1]\),
	\[
	\frac14\,rN_A(r)
	\leq
	\bigl|
	\{x\in[0,1]:
	\operatorname{dist}(x,A)<r\}
	\bigr|
	\leq
	3rN_A(r).
	\]
\end{lemma}
\begin{proof}
	Let \(I_1,\ldots,I_{N_A(r)}\) be intervals of length \(r\) covering
	\(A\). Enlarging each \(I_i\) by \(r\) in both directions gives an
	interval of length \(3r\), and the enlarged intervals cover the
	\(r\)-neighborhood of \(A\). This proves the upper bound.
	
	For the lower bound, choose a maximal \(2r\)-separated set
	\(\{x_1,\ldots,x_m\}\subset A\). Then the sets
	\((x_i-r,x_i+r)\cap[0,1]\) are pairwise disjoint, contained in the
	\(r\)-neighborhood of \(A\), and have length at least \(r\). Hence
	the \(r\)-neighborhood has measure at least \(mr\).
	
	By maximality, the intervals of radius \(2r\) centered at the \(x_i\)
	cover \(A\), so \(N_A(4r)\leq m\). Since an interval of length \(4r\)
	can be covered by four intervals of length \(r\),
	\(N_A(r)\leq4N_A(4r)\leq4m\). Thus \(m\geq N_A(r)/4\), which proves
	the lower bound.
\end{proof}

Since \(A_{\cL}\) is countable and hence has Lebesgue measure zero,
\[
    V_{\cL}(r)
    =
    \bigl|\{x\in[0,1]:
    \operatorname{dist}(x,A_{\cL})<r\}\bigr|.
\]
Applying Lemma~\ref{lem:tube-cover} to the set \(A=A_{\cL}\) induced by a normalized fractal string $\cL$ gives
\begin{equation}\label{eq:tube-cover}
	V_{\cL}(r)\asymp rN_{A_{\cL}}(r),
	\qquad r\downarrow0.
\end{equation}
In fact, the two constants implicit in
\eqref{eq:tube-cover} may be chosen independently of the string
\(\cL\).

The lower and upper Minkowski (or box counting) dimensions of a compact set
\(A\subset\R\) are defined by
\[
\underline{\dim}_M A
:=
\liminf_{r\downarrow0}
\frac{\log N_A(r)}{\log(1/r)}\qquad\text{and}\qquad
\overline{\dim}_M A
:=
\limsup_{r\downarrow0}
\frac{\log N_A(r)}{\log(1/r)}.
\]
If these two quantities agree, their common value is denoted by
\(\dim_M A\) and is called the Minkowski dimension of \(A\),
see \cite[Section~3.1]{Falconer2003} or \cite[Section~1.1]{LapidusVanFrankenhuijsen2000}. By \eqref{eq:tube-cover}, the same dimensions can be read from the tube function (which is a well-known fact).
In particular, if \(\dim_M A_{\cL}=D\) exists, then $V_{\cL}(r)=r^{1-D+o(1)}$ as $r\downarrow0$. In this context, the Minkowski dimension of $A_\cL$ is also referred to as the Minkowski dimension of the fractal string $\cL$.

\section{Main Results}\label{sec:Results}

\subsection{The area deficit and the integrated tube volume}

Given a normalized fractal string \(\cL\), let
\(a=a(\cL)\) be its tail-sum sequence as defined in \eqref{eq:tails}, and set
\[
K(\cL):=K(a(\cL)),
\]
where $K(a(\cL))$ is defined by \eqref{eq:Ka}. To describe the area deficit between $K(\cL)$ and its floating body we use the function
\begin{equation}\label{eq:Fdef}
    \cF_{\cL}(\eps)
    :=\eps^3\int_\eps^1\frac{V_{\cL}(r)}{r^2}\,dr,
    \qquad 0<\eps<1.
\end{equation}
In view of \eqref{eq:tube-cover}, the function \(\cF_{\cL}\) may
equivalently be viewed, up to multiplicative constants, as
$\eps^3\int_\eps^1 N_{A_{\cL}}(r)\,\frac{dr}{r}$.
Since \(dr/r\) is the Haar measure on the multiplicative group
\((0,\infty)\), this expression records the covering complexity of
\(A_{\cL}\) integrated over logarithmic scales.
We remark that the choice of the upper integration bound is somewhat
arbitrary. Replacing \(1\) by any fixed positive number changes
\(\cF_{\cL}(\eps)\) by at most \(O(\eps^3)\), which is negligible for
the asymptotic regimes considered below.

\begin{theorem}\label{thm:main}
Let $\cL$ be an infinite normalized fractal string, and let $K(\cL)$ be its associated almost polygonal body.  Then
\[
\Delta_{K(\cL)}(\eps^3)
\asymp
\cF_{\cL}(\eps),
\qquad \eps\downarrow0.
\]
\end{theorem}
\begin{proof}
Put \(\delta=\eps^3\) and let \(m=n(\eps^3)\), where \(n(\cdot)\)
is defined by \eqref{eq:ndef}. Since \(m\to\infty\) as
\(\eps\downarrow0\), we may assume that \(m\geq2\). Then $\ell_m\le 2^{1/3}\eps<\ell_{m-1}$. For $2\le k\le m$, set
\[
    Q_k(\eps)
    :=\frac{\ell_{k-1}\ell_{k-2}(\ell_{k-1}+\ell_{k-2})}{2\eps^3}.
\]
Since $\ell_{k-2}\ge\ell_{k-1}>2^{1/3}\eps$, we have $Q_k(\eps)\ge2$.
Moreover, since
\[
    \log\frac{\ell_{k-1}}\eps
    \leq
    \log\frac{\ell_{k-1}+\ell_{k-2}}{2\eps}
    \leq
    \log\frac{\ell_{k-2}}\eps
\]
and both endpoint expressions are at least
\(\frac13\log 2\), we obtain
\[
    \frac14\left(
        1+\log\frac{\ell_{k-1}}\eps
        +\log\frac{\ell_{k-2}}\eps
    \right)
    \leq
    \log Q_k(\eps)
    \leq
    2\left(
        1+\log\frac{\ell_{k-1}}\eps
        +\log\frac{\ell_{k-2}}\eps
    \right).
\]
Summing over \(k=2,\ldots,m\), and observing that each term
\(\log(\ell_j/\eps)\), \(0\leq j<m\), occurs once or twice
among the first two logarithmic summands in the expansion of
\(\log Q_k(\eps)\), we conclude that
\[
    \sum_{k=2}^m \log Q_k(\eps)
    \asymp
    \sum_{j=0}^{m-1}
    \left(1+\log\frac{\ell_j}{\eps}\right),
    \qquad \eps\downarrow0.
\]
Multiplying both sides by $\eps^3$ provides an asymptotic reformulation of the third term in the upper bound of Proposition~\ref{prop:SW}. In particular, it is now transparent that the second term
\(\delta n(\delta)=\eps^3m\) in the upper bound is absorbed by the third term.  Since by \eqref{eq:tails} we have $ a_m=\sum_{j=m}^\infty\ell_j$,
Proposition~\ref{prop:SW} gives
\begin{equation}\label{eq:discrete-profile}
    \Delta_{K(\cL)}(\eps^3)
    \asymp
    \eps^2\sum_{j=m}^\infty\ell_j
    +\eps^3\sum_{j=0}^{m-1}
       \left(1+\log\frac{\ell_j}{\eps}\right), \qquad\eps\downarrow 0.
\end{equation}
To complete the proof, it remains to identify the right-hand side with \eqref{eq:Fdef}.  By Tonelli's theorem and \eqref{eq:tube-sum},
\[
    \cF_{\cL}(\eps)
    =\sum_{j=0}^\infty
      \eps^3\int_\eps^1\frac{\min\{2r,\ell_j\}}{r^2}\,dr.
\]
For a fixed length $\ell\le1$, elementary integration gives
\begin{equation*}
 \eps^3\int_\eps^1\frac{\min\{2r,\ell\}}{r^2}\,dr
 =
 \begin{cases}
   \eps^2\ell(1-\eps),&\ell\le2\eps,\\[1ex]
   \eps^3\bigl(2\log(\ell/(2\eps))+2-\ell\bigr),&\ell>2\eps.
 \end{cases}
\end{equation*}
For \(0<\eps<1/2\), we have
\[
    \frac12\ell
    \leq \ell(1-\eps)\leq\ell.
\]
Moreover, whenever \(2\eps<\ell\leq1\),
\[
    \frac12\left(1+\log\frac{\ell}{\eps}\right)
    \leq
    2\log\frac{\ell}{2\eps}+2-\ell
    \leq
    2\left(1+\log\frac{\ell}{\eps}\right).
\]
Consequently,
\begin{equation}\label{eq:Fdiscrete}
    \cF_{\cL}(\eps)
    \asymp
    \eps^2\sum_{\ell_j\le2\eps}\ell_j
    +\eps^3\sum_{\ell_j>2\eps}
       \left(1+\log\frac{\ell_j}\eps\right),\qquad\eps\downarrow 0.
\end{equation}
For lengths in the intermediate range $2^{1/3}\eps<\ell\leq2\eps$,
the two possible contributions satisfy
\[
\eps^2\ell\asymp\eps^3
\asymp
\eps^3\left(1+\log\frac{\ell}{\eps}\right).
\]
Thus moving the threshold between \(2^{1/3}\eps\) and \(2\eps\)
changes the resulting expression only by absolute multiplicative
constants.  So, comparing \eqref{eq:discrete-profile} and \eqref{eq:Fdiscrete} proves the theorem.
\end{proof}

\begin{remark}\label{rem:finite}
\begin{itemize}
    \item[(i)] If one also allows finite normalized strings
	\(\cL=(\ell_0,\ldots,\ell_{N-1})\), then $V_{\cL}(r)=2Nr$ for all sufficiently small $r$, and $\cF_{\cL}(\eps)\asymp \eps^3\log(1/\eps)$.  Since $\delta=\varepsilon^3$, this is the same order $\delta\log(1/\delta)$ known for the area deficit of planar polygons from \cite{Schuett1991,SchuettWerner1990}.

    \item[(ii)] For a convex body \(K\subset\R^2\), let \(K^\delta\) denote its illumination
body, consisting of all points \(x\in\R^2\) for which
\(\vol(\operatorname{conv}(K\cup\{x\})\setminus K)\leq\delta\).
Estimates for the area difference
\(\vol(K^\delta\setminus K)\) in the case that $K$ is an almost polygonal body
were obtained in \cite[Corollary~3]{Werner1997}.
Writing \(m=\min\{n:\ell_n\leq\eps\}\), these estimates give
\[
    \vol(K(\cL)^{\eps^3}\setminus K(\cL))
    \asymp
    \eps^2a_m+\eps^3m
    =
    \eps^2V_{\cL}(\eps/2)
    \asymp
    \eps^2V_{\cL}(\eps).
\]
This relation allows, in principle, to derive analogous dimension estimates and to transfer the examples below to illumination bodies, although the logarithmic and slowly varying corrections may differ, since the illumination-body area deficit depends on the tube function at a single scale rather than on its integral over a range of scales as encoded in $f_{\cL}$.
\end{itemize}
\end{remark}

Theorem~\ref{thm:main} immediately allows us to identify the power exponent of the area deficit for almost polygonal bodies.

\begin{corollary}\label{cor:dimension}
Let $\cL$ be an infinite normalized fractal string, and let $K(\cL)$ be the associated almost polygonal body. Let $\underline D$ and $\overline D$ be the lower and upper Minkowski dimensions of $A_{\cL}$.  Then
\begin{align}\label{eq:cor1eq1}
 \liminf_{\delta\downarrow0}
 \frac{\log\Delta_{K(\cL)}(\delta)}{\log\delta}
 =1-\frac{\overline D}{3}\qquad\text{and}\qquad
 \limsup_{\delta\downarrow0}
 \frac{\log\Delta_{K(\cL)}(\delta)}{\log\delta}
 =1-\frac{\underline D}{3}.
\end{align}
In particular, if $D=\dim_M A_{\cL}$ exists, then $\Delta_{K(\cL)}(\delta)=\delta^{1-D/3+o(1)}$ as $\delta\downarrow 0$.
\end{corollary}
\begin{proof}
	Put
	\[
	I_{\cL}(\eps):=\int_\eps^1 N_{A_{\cL}}(r)\,\frac{dr}{r}.
	\]
	Since \(N_{A_{\cL}}\) is nonincreasing,
	\[
	I_{\cL}(\eps)
	\geq
	\int_\eps^{2\eps}
	N_{A_{\cL}}(r)\,\frac{dr}{r}
	\geq
	(\log2)N_{A_{\cL}}(2\eps).
	\]
	Every interval of length \(2\eps\) can be covered by two intervals of
	length \(\eps\), so
	\(N_{A_{\cL}}(\eps)\leq2N_{A_{\cL}}(2\eps)\). Hence
	\[
	\frac{\log2}{2}N_{A_{\cL}}(\eps)
	\leq
	I_{\cL}(\eps)
	\leq
	N_{A_{\cL}}(\eps)\log(1/\eps).
	\]
	Taking logarithms gives
\[
    \log N_{A_{\cL}}(\eps)+\log\frac{\log 2}{2}
    \leq
    \log I_{\cL}(\eps)
    \leq
    \log N_{A_{\cL}}(\eps)+\log\log(1/\eps).
\]
	Since
	\(\log\log(1/\eps)/\log(1/\eps)\to0\),
	\(I_{\cL}\) and \(N_{A_{\cL}}\) have the same lower and upper
	logarithmic orders.
	
	By Theorem~\ref{thm:main},
	\(\Delta_{K(\cL)}(\eps^3)\asymp\cF_{\cL}(\eps)\).
	Moreover, \eqref{eq:Fdef} and \eqref{eq:tube-cover} give
	\(\cF_{\cL}(\eps)\asymp\eps^3I_{\cL}(\eps)\). Therefore $\Delta_{K(\cL)}(\eps^3)\asymp\eps^3I_{\cL}(\eps)$.
	Writing \(\delta=\eps^3\), we obtain
	\[
	\frac{\log\Delta_{K(\cL)}(\delta)}{\log\delta}
	=
	1-\frac13
	\frac{\log I_{\cL}(\eps)}{\log(1/\eps)}
	+o(1).
	\]
	Taking lower and upper limits yields \eqref{eq:cor1eq1}.
\end{proof}

\begin{remark}\label{rem:strict-interpolation}
The area deficits of almost polygonal bodies associated with infinite
fractal strings lie strictly between the polygonal and smooth orders.
Indeed, by \eqref{eq:tube-sum},
\[
    V_{\cL}(r)\longrightarrow0,
    \qquad
    \frac{V_{\cL}(r)}r\longrightarrow\infty,
    \qquad r\downarrow0.
\]
The first limit follows by dominated convergence, whereas the
second follows from the fact that the string contains infinitely
many positive lengths. Consequently,
\[
    \frac{\cF_{\cL}(\eps)}
         {\eps^3\log(1/\eps)}
    =
    \frac1{\log(1/\eps)}
    \int_\eps^1\frac{V_{\cL}(r)}r\,\frac{dr}{r}
    \longrightarrow\infty,
\]
since the logarithmic average of a function tending to infinity
also diverges. On the other hand, the substitution \(s=\eps/r\)
gives
\[
    \frac{\cF_{\cL}(\eps)}{\eps^2}
    =
    \int_\eps^1 V_{\cL}(\eps/s)\,ds
    \longrightarrow0
\]
by dominated convergence. Together with
Theorem~\ref{thm:main}, these relations yield
\[
    \delta\log(1/\delta)
    =
    o\bigl(\Delta_{K(\cL)}(\delta)\bigr)
    \qquad\text{and}\qquad
    \Delta_{K(\cL)}(\delta)
    =
    o(\delta^{2/3}).
\]
\end{remark}

The Minkowski dimension of \(A_{\cL}\) alone does not determine
logarithmic or other slowly varying factors of the area deficit. The
natural refinement is a tube estimate for the fractal string. For this purpose we recall that a positive measurable and locally integrable function \(L:[1,\infty)\to(0,\infty)\) is slowly varying at infinity if
\(L(cx)/L(x)\to1\) as \(x\to\infty\) for every \(c>0\), see e.g.~\cite{BinghamGoldieTeugels}.

\begin{corollary}\label{cor:regular}
Let $\cL$ be an infinite normalized fractal string, and let $K(\cL)$ be the associated almost polygonal body. Assume that, for some \(D\in[0,1]\),
\begin{equation}\label{eq:Vregular}
    V_{\cL}(r)\asymp r^{1-D}L(1/r),\qquad r\downarrow 0,
\end{equation}
where $L$ is slowly varying.
\begin{itemize}
\item[(i)] If $0<D\le1$, then $\Delta_{K(\cL)}(\delta)\asymp\delta^{1-D/3}L(\delta^{-1/3})$.
\item[(ii)] If $D=0$, then $\Delta_{K(\cL)}(\delta)\asymp\delta\,\widehat L(\delta^{-1/3})$, where  $\widehat L(x):=\int_1^xL(u)\,\frac{du}{u}$.
\end{itemize}
\end{corollary}
\begin{proof}
	The assumption \eqref{eq:Vregular} means that there exist  \(r_0\in(0,1)\) and positive constants $c,C$ such that
    $$
    c r^{1-D} L(1/r)\leq V_{\cL}(r)\leq C r^{1-D} L(1/r)
    $$
    holds for
	\(0<r\leq r_0\). Multiplying by $r^{-2}$ and integrating over the interval $(\eps,r_0]$ for $\eps>0$ shows that the contribution of
	\((\eps,r_0]\) to \(\cF_{\cL}(\eps)\) is comparable to
	\[
	\eps^3\int_\eps^{r_0}r^{-1-D}L(1/r)\,dr,
	\]
	while the contribution of \([r_0,1]\) is \(O(\eps^3)\).
	
	Suppose first that \(D>0\). After the substitution \(u=1/r\),
	Karamata's theorem
	\cite[Theorem~1.5.11(i), p.~28]{BinghamGoldieTeugels} gives
\[
\int_\eps^{r_0}r^{-1-D}L(1/r)\,dr
\asymp
\eps^{-D}L(1/\eps).
\]
	Since the right-hand side tends to infinity, the  contribution of	\((\eps,r_0]\) to \(\cF_{\cL}(\eps)\) dominates that of \([r_0,1]\). Hence
	\(\cF_{\cL}(\eps)\asymp\eps^{3-D}L(1/\eps)\), as $\eps\downarrow 0$. Theorem~\ref{thm:main},
	with \(\eps=\delta^{1/3}\), proves part~\textup{(i)}.
	
	Now let \(D=0\). By \eqref{eq:Vregular} and \eqref{eq:tube-cover},
	\(L(1/r)\asymp N_{A_{\cL}}(r)\). Since \(\cL\) is infinite,
	\(N_{A_{\cL}}(r)\to\infty\) as \(r\downarrow0\), and therefore
	\(L(x)\to\infty\) and \(\widehat L(x)\to\infty\) as $x\to\infty$. Thus
	\[
	\int_\eps^{r_0}\frac{L(1/r)}r\,dr
	=
	\widehat L(1/\eps)-\widehat L(1/r_0)
	\asymp
	\widehat L(1/\eps).
	\]
	Consequently,
	\(\cF_{\cL}(\eps)\asymp\eps^3\widehat L(1/\eps)\), and another
	application of Theorem~\ref{thm:main} proves part~\textup{(ii)}.
\end{proof}

\subsection{A realization criterion}
\label{subsec:realization}

The preceding results express the floating-body area deficit of an almost
polygonal body in terms of the tube geometry of its associated fractal
string. We now address the converse problem of realizing prescribed
asymptotic deficit rates. Tube asymptotics of fractal strings, including
generalized Minkowski contents and regularly varying gauge functions,
have been studied extensively, see, for example,
\cite{EichingerWinter2019,HeLapidus1997,LapidusVanFrankenhuijsen2000,RatajWinter2013}.
For our purposes, the following elementary approximation result is
sufficient.

\begin{lemma}\label{lem:tube-realization}
	Let \(v:[0,r_0]\to[0,\infty)\) be continuous, nondecreasing and
	concave, with \(v(0)=0\) and $v(r)/r\to\infty$ as $r\downarrow0$.
	Then there is an infinite normalized fractal string \(\cL\) such that
	\[
	V_{\cL}(r)=v(r)+O(r),
	\qquad r\downarrow0.
	\]
\end{lemma}
\begin{proof}
	After decreasing \(r_0\), we may assume that \(v(r_0)<1\). Since \(v\) is
	concave, its right derivative \(v'_+\) is nonnegative and
	nonincreasing. For \(0<s<2r_0\), put
	\[
	M(s):=
	\left\lfloor\frac12v'_+(s/2)\right\rfloor
	\]
	and let $\lambda_j:=\bigl|\{s\in(0,2r_0):M(s)\geq j\}\bigr|$, for $j\geq 1$.
	Since \(v(r)/r\to\infty\), concavity implies
	\(v'_+(r)\to\infty\) as \(r\downarrow0\). Hence
	\(\lambda_j>0\) for every \(j\geq1\), and the sequence
	\((\lambda_j)_{j\geq1}\) is nonincreasing. Moreover, by Tonelli's theorem,
	\[
	S:=\sum_{j=1}^{\infty}\lambda_j
	=
	\int_0^{2r_0}M(s)\,ds
	\leq v(r_0)<1.
	\]
	Adjoin the length \(1-S\) to the sequence
	\((\lambda_j)_{j\geq1}\) and rearrange the resulting sequence in
	nonincreasing order. Since
	\[
	(1-S)+\sum_{j=1}^{\infty}\lambda_j=1,
	\]
	this defines an infinite normalized fractal string \(\cL\). For all sufficiently small
	\(r\), the added interval contributes \(2r\), whereas the remaining
	intervals contribute \(\int_0^{2r}M(s)\,ds\). Since
	\[
	\frac12v'_+(s/2)-1
	\leq M(s)
	\leq\frac12v'_+(s/2),
	\]
	we obtain $0\leq V_{\cL}(r)-v(r)\leq2r$.
	The assumption \(v(r)/r\to\infty\) as $r\downarrow0$ also implies that the constructed
	string is infinite. This concludes the proof.
\end{proof}

Let \(\delta_0>0\), and let
\(h:(0,\delta_0)\to(0,\infty)\) be a positive \(C^1\)-function.
We want to understand when \(h\) can occur, up to multiplicative constants, as the deficit rate of some (almost polygonal) convex body generated by some fractal string.
In view of \eqref{eq:Fdef} and Theorem~\ref{thm:main}, this amounts to asking which volume function
\(V_{\cL}\) would produce the prescribed quantity \(h(\eps^3)\) after
integration. This suggests differentiating the normalized function
\(h(t^3)/t^3\).

Accordingly, for \(0<t<\delta_0^{1/3}\), define
\begin{equation}\label{eq:target-tube}
	v_h(t)
	:=
	-t^2\frac{d}{dt}
	\left(\frac{h(t^3)}{t^3}\right)
	=
	3\frac{h(t^3)}{t^2}-3t h'(t^3).
\end{equation}
Indeed, if for some fractal string \(\cL\) one had the exact identity
\[
h(\eps^3)=\cF_{\cL}(\eps)
=
\eps^3\int_\eps^1\frac{V_{\cL}(r)}{r^2}\,dr,
\]
then differentiation would give
\(v_h(\eps)=V_{\cL}(\eps)\). Thus \(v_h\) is the function that is
naturally associated with a prescribed deficit rate \(h\) when the
transform in \eqref{eq:Fdef} is inverted.

\begin{theorem}\label{thm:realization}
	Let \(\delta_0>0\), and let \(h:(0,\delta_0)\to(0,\infty)\) be of class \(C^1\). Suppose that,
	for all sufficiently small \(t>0\), the function \(v_h\) in
	\eqref{eq:target-tube} is nondecreasing and concave, and that
	\[
	v_h(t)\longrightarrow0,
	\qquad
	\frac{v_h(t)}t\longrightarrow\infty,
	\qquad t\downarrow0.
	\]
	Then there is an infinite normalized fractal string \(\cL_h\) such that
	\[
	\cF_{\cL_h}(\eps)
	=
	h(\eps^3)(1+o(1))
	\]
	and consequently $\Delta_{K(\cL_h)}(\delta)\asymp h(\delta)$.
	If
	\(\log h(\delta)/\log\delta\to\alpha\) as $\delta\downarrow0$, then
	\(\dim_M A_{\cL_h}=3(1-\alpha)\).
\end{theorem}

\begin{proof}
	Choose \(t_0>0\) sufficiently small such that \(v_h\) is
	nondecreasing and concave on \((0,t_0]\), and extend it continuously
	to zero by setting \(v_h(0)=0\). Lemma~\ref{lem:tube-realization}
	then yields an infinite normalized fractal string \(\cL_h\) satisfying
	\[
	V_{\cL_h}(r)=v_h(r)+O(r),
	\qquad r\downarrow0.
	\]
	By \eqref{eq:target-tube},
	\[
	\frac{h(\eps^3)}{\eps^3}
	=
	\frac{h(t_0^3)}{t_0^3}
	+
	\int_\eps^{t_0}\frac{v_h(r)}{r^2}\,dr.
	\]
	Since \(v_h(r)/r\to\infty\), for every \(M>0\) we have
	\(v_h(r)\geq Mr\) for all sufficiently small \(r\). It follows that
	\[
	\frac{1}{\log(1/\eps)}
	\int_\eps^{t_0}\frac{v_h(r)}{r^2}\,dr
	\longrightarrow\infty,\qquad\eps\downarrow0.
	\]
	In particular, $\eps^3\log(1/\eps)=o(h(\eps^3))$.
	
	The \(O(r)\)-error in the approximation of \(V_{\cL_h}(r)\) therefore
	contributes only
	\(O(\eps^3\log(1/\eps))=o(h(\eps^3))\) to
	\(\cF_{\cL_h}(\eps)\). The part of the integral in
	\eqref{eq:Fdef} over the fixed interval \([t_0,1]\) contributes
	\(O(\eps^3)\), which is negligible as well. Consequently,
	\[
	\cF_{\cL_h}(\eps)
	=
	h(\eps^3)(1+o(1)).
	\]
	Theorem~\ref{thm:main} gives
	\(\Delta_{K(\cL_h)}(\delta)\asymp h(\delta)\). If the logarithmic
	exponent of \(h\) exists, Corollary~\ref{cor:dimension} then yields
	\(\dim_M A_{\cL_h}=3(1-\alpha)\).
\end{proof}

\begin{remark}\label{rem:realization-comparison}
	The criterion is sharp at the level of the integrated tube functions.
	Indeed, for an infinite normalized fractal string \(\cL\), define
	\(h_{\cL}(\delta):=\cF_{\cL}(\delta^{1/3})\). Then
	\(h_{\cL}\) is of class \(C^1\) and satisfies
	\(v_{h_{\cL}}=V_{\cL}\). By \eqref{eq:tube-sum},
	\(V_{\cL}\) is nondecreasing and concave, tends to zero, and satisfies
	\(V_{\cL}(r)/r\to\infty\). Thus, up to replacing a deficit rate by a
	comparable \(C^1\)-function, Theorem~\ref{thm:realization}
	characterizes the asymptotic orders generated by infinite strings.
	
	The relation with \cite[Theorem~2.2]{SchuettWerner1992} is particularly
	simple. With
	\[
	f(t)=3h'(t^3)-3\frac{h(t^3)}{t^3}
	\]
	as in that theorem, \eqref{eq:target-tube} reads \(v_h(t)=-tf(t)\).
	Hence the assumption in \cite[Theorem~2.2]{SchuettWerner1992} that
	\(tf(t)\) is decreasing and strictly convex corresponds to
	\(v_h\) being nondecreasing and strictly concave. In the almost
	polygonal branch of the proof of \cite[Theorem~2.2]{SchuettWerner1992},
	the additional assumption
	\(f(t)/\log t\to c\in(0,\infty]\) implies
	\(v_h(t)/t\to\infty\). Theorem~\ref{thm:realization}, by contrast,
	requires only \(v_h(t)/t\to\infty\) together with non-strict
	concavity.
	
	For instance, if
	\(h(t)=t(\log(1/t))^q\) with \(1<q<2\), then $v_h(t)=3^q q\,t(\log(1/t))^{q-1}$.
	Hence the assumptions of Theorem~\ref{thm:realization} hold for all
	sufficiently small \(t\), whereas
	\(f(t)/\log t\to0\). Thus this rate is not covered by the hypotheses of
	\cite[Theorem~2.2]{SchuettWerner1992}. It also follows directly from
	Example~\ref{ex:stretch} below.
\end{remark}

\section{Examples}\label{sec:Ex}

In this section we illustrate our results through several examples,
recovering in particular the three cases considered in
\cite{SchuettWerner1992}.

\begin{example}[Polynomial strings]\label{ex:power}
	Let
	\[
	a_n=(n+1)^{-\beta},\qquad \beta>0.
	\]
	Then \(\ell_n=a_n-a_{n+1}\asymp n^{-\beta-1}\), so the number of
	lengths exceeding \(r\) is comparable to \(r^{-1/(\beta+1)}\).
	Moreover, the sum of the remaining lengths is of order
	\(r^{\beta/(\beta+1)}\). Hence the tube formula gives
	\(V_{\cL}(r)\asymp r^{\beta/(\beta+1)}\). Therefore
	\(D:=\dim_M A_{\cL}=1/(\beta+1)\), and
	\[
	\Delta_{K(\cL)}(\delta)
	\asymp
	\delta^{1-D/3}
	=
	\delta^{(2+3\beta)/(3+3\beta)}.
	\]
	This recovers the exponent in Example~(2.8) of
	\cite{SchuettWerner1992}.
\end{example}

\begin{example}[Stretched-exponential strings]\label{ex:stretch}
Fix $p>0$ and, after normalization, let
\[
    \ell_n\asymp e^{-n^p}.
\]
The number of lengths exceeding $r$ is comparable to $(\log(1/r))^{1/p}$.  The tail of the length sequence is of smaller order than $r(\log(1/r))^{1/p}$, and hence $V_{\cL}(r)\asymp r(\log(1/r))^{1/p}$.
Thus \(\dim_M A_{\cL}=0\), and  Corollary~\ref{cor:regular}(ii) yields
\begin{equation*}
    \Delta_{K(\cL)}(\delta)
    \asymp
    \delta\bigl(\log(1/\delta)\bigr)^{1+1/p}.
\end{equation*}
For \(p=1\), the gaps decay exponentially and one recovers
Example~(2.9) of \cite{SchuettWerner1992}, while varying \(p>0\)
produces a continuous family of dimension-zero logarithmic
corrections.
\end{example}

\begin{example}[Logarithmically corrected power strings]\label{ex:logpower}
Let $0<D<1$ and $\gamma\in\R$, and suppose that, after a finite adjustment and normalization,
\[
    \ell_n\asymp n^{-1/D}(\log n)^{-\gamma/D}.
\]
The number of lengths exceeding \(r\) is comparable to
\(r^{-D}(\log(1/r))^{-\gamma}\). By the integral test, the sum
of the remaining lengths has the same order as \(r\) times this
number. Hence the tube formula gives
\(V_{\cL}(r)\asymp r^{1-D}(\log(1/r))^{-\gamma}\). Consequently,
\begin{equation*}
    \Delta_{K(\cL)}(\delta)
    \asymp
    \delta^{1-D/3}(\log(1/\delta))^{-\gamma}.
\end{equation*}
This family realizes arbitrary logarithmic corrections at every dimension $D\in(0,1)$.
\end{example}

\begin{example}[Iterated-logarithmic strings]
	\label{ex:iterated-log}
	For \(j\geq1\), define the iterated logarithm by
	\[
	\log_1 x:=\log x,
	\qquad
	\log_{j+1}x:=\log(\log_jx).
	\]
	Let \(m\geq1\) and \(\gamma>0\), and suppose that, after a finite
	adjustment and normalization,
	\[
	\ell_n
	\asymp
	\frac{1}{
		n\bigl(\prod_{j=1}^{m-1}\log_j n\bigr)
		(\log_m n)^{1+\gamma}
	},
	\]
	where the empty product is understood as \(1\). The integral test gives $\sum_{k>n}\ell_k\asymp(\log_m n)^{-\gamma}$,
	and this term dominates \(n\ell_n\). Consequently, $V_{\cL}(r)\asymp\bigl(\log_m(1/r)\bigr)^{-\gamma}$,
	so that \(\dim_M A_{\cL}=1\) and
	\[
	\Delta_{K(\cL)}(\delta)
	\asymp
	\frac{\delta^{2/3}}
	{(\log_m(1/\delta))^\gamma}.
	\]
	For \(m=2\) and \(\gamma=1\), the tail sums satisfy $a_n\asymp 1/\log\log n$,
	which recovers Example~(2.10) of \cite{SchuettWerner1992}, up to a
	finite adjustment.
\end{example}

\begin{example}[A string without Minkowski dimension]
	\label{ex:no-minkowski-dimension}
	
	Let \(Z:=\sum_{k=2}^{\infty}2^{-k!/2}\). For every \(k\geq2\), let
	the length $r_k:=Z^{-1}2^{-k!}$
	occur with multiplicity \(m_k:=2^{k!/2}\). Arranging these lengths
	in nonincreasing order defines a normalized fractal string \(\cL\),
	since \(\sum_{k=2}^{\infty}m_kr_k=1\).
	
	We claim that
	\(\underline{\dim}_M A_{\cL}=0\) and
	\(\overline{\dim}_M A_{\cL}=1/2\). Put
	\(M_k:=\sum_{j=2}^k m_j\) and
	\(R_{k+1}:=\sum_{j=k+1}^{\infty}m_jr_j\).
	The factorial separation of the scales gives
	\(M_k\asymp r_k^{-1/2}\) and
	\(R_{k+1}\asymp r_{k+1}^{1/2}\).
	
	For \(r_{k+1}\leq r\leq r_k\), the blocks with index \(j\leq k\)
	contribute a quantity comparable to \(rM_k\) to the tube volume,
	whereas the remaining blocks contribute \(R_{k+1}\). Hence
	\eqref{eq:tube-cover} gives
	\begin{equation}\label{eq:lacunary-covering}
		N_{A_{\cL}}(r)
		\asymp
		r_k^{-1/2}+\frac{r_{k+1}^{1/2}}{r},
		\qquad r_{k+1}\leq r\leq r_k.
	\end{equation}
	Taking \(r=r_k\) shows that
	\(N_{A_{\cL}}(r_k)\asymp r_k^{-1/2}\). On the other hand,
	\eqref{eq:lacunary-covering} gives
	\(N_{A_{\cL}}(r)=O(r^{-1/2})\), and therefore
	\(\overline{\dim}_M A_{\cL}=1/2\).
	
	To determine the lower Minkowski dimension, put
	\(s_k:=r_k^{3/4}r_{k+1}^{1/4}\). Then
	\(r_{k+1}<s_k<r_k\), and the second term in
	\eqref{eq:lacunary-covering} is of smaller order than the first.
	Thus \(N_{A_{\cL}}(s_k)\asymp r_k^{-1/2}\). Since
	\[
	\log\frac1{r_k}=k!\log2+O(1),
	\qquad
	\log\frac1{s_k}
	=
	\frac34\log\frac1{r_k}
	+\frac14\log\frac1{r_{k+1}},
	\]
	we obtain
	\(\log N_{A_{\cL}}(s_k)/\log(1/s_k)\to0\), and hence
	\(\underline{\dim}_M A_{\cL}=0\).
	
	Consequently, Corollary~\ref{cor:dimension} yields
	\[
	\liminf_{\delta\downarrow0}
	\frac{\log\Delta_{K(\cL)}(\delta)}{\log\delta}
	= 1 - {\frac{1/2}{3}} =
	\frac56,
	\qquad
	\limsup_{\delta\downarrow0}
	\frac{\log\Delta_{K(\cL)}(\delta)}{\log\delta}
	= 1- {\frac 0 3} =
	1.
	\]
In particular, there is no exponent \(\alpha\) such that $\Delta_{K(\cL)}(\delta)=\delta^{\alpha+o(1)}$ as $\delta\downarrow0$.
\end{example}

\section{Consequences for random polygons}
\label{sec:random-polygons}

Let \(\cL\) be an infinite normalized fractal string and let
\(K(\cL)\) be the associated almost polygonal body. Let
\(X_1,X_2,\ldots\) be independent random points uniformly distributed
in \(K(\cL)\), and define the random polygon
\[
P_n:=\operatorname{conv}(X_1,\ldots,X_n),\qquad n\geq 3.
\]
We are interested in the expected missed area $\mathbb E\vol\bigl(K(\cL)\setminus P_n\bigr)$ and the expected number of vertices \(\mathbb{E}f_0(P_n)\) of  \(P_n\) as $n\to\infty$. We refer to \cite{Reitzner2010} for general background on random polytopes and their connection with floating bodies.

\begin{corollary}\label{cor:random-polygons}
	Let $\cL$ be an infinite normalized fractal string, and let $K(\cL)$ be its associated almost polygonal body.
	Let \(\underline D\) and \(\overline D\) denote the lower and upper
	Minkowski dimensions of \(A_{\cL}\). Then
	\[
	\liminf_{n\to\infty}
	\frac{
		\log\mathbb E\vol(K(\cL)\setminus P_n)
	}{\log 1/n}
	=
	1-\frac{\overline D}{3}\qquad\text{and}\qquad
	\limsup_{n\to\infty}
	\frac{
		\log\mathbb E\vol(K(\cL)\setminus P_n)
	}{\log 1/n}
	=
	1-\frac{\underline D}{3}.
	\]
	Moreover,
	\[
	\liminf_{n\to\infty}
	\frac{\log\mathbb E f_0(P_n)}{\log n}
	=
	\frac{\underline D}{3}\qquad\text{and}\qquad
	\limsup_{n\to\infty}
	\frac{\log\mathbb E f_0(P_n)}{\log n}
	=
	\frac{\overline D}{3}.
	\]
	In particular, if \(D=\dim_M A_{\cL}\) exists, then
	\[
	\mathbb E\vol\bigl(K(\cL)\setminus P_n\bigr)
	=
	n^{-1+D/3+o(1)}\qquad\text{and}\qquad
	\mathbb E f_0(P_n)
	=
	n^{D/3+o(1)}.
	\]
\end{corollary}
\begin{proof}
	Set
	\[
	v:=\vol(K(\cL)),
	\qquad
	\delta_n:=\frac{v}{n},
	\qquad
	M_n:=\mathbb E\vol\bigl(K(\cL)\setminus P_n\bigr).
	\]
	Theorem~1 of  \cite{BaranyLarman1988} gives
	\begin{equation}\label{eq:random-wet-part}
		M_n
		\asymp
		\Delta_{K(\cL)}(\delta_n),
		\qquad n\to\infty,
	\end{equation}
	see also \cite[page 321]{SchneiderWeil}.
	
	Note that the lower and upper logarithmic orders of
	\(\Delta_{K(\cL)}\) may be computed along the sequence
	\((\delta_n)_{n\geq1}\). Indeed, this follows since
	\(\delta_{n+1}/\delta_n\to1\) as \(n\to\infty\) and since
	\(\delta\mapsto\Delta_{K(\cL)}(\delta)\) is nondecreasing, meaning that
\(\delta_{n+1}\leq\delta\leq\delta_n\) implies
	\[
	\Delta_{K(\cL)}(\delta_{n+1})
	\leq
	\Delta_{K(\cL)}(\delta)
	\leq
	\Delta_{K(\cL)}(\delta_n).
	\]
	It follows from Corollary~\ref{cor:dimension} that
	\begin{align*}
		\liminf_{n\to\infty}
		\frac{\log\Delta_{K(\cL)}(\delta_n)}
		{\log\delta_n}
		&=
		\liminf_{\delta\downarrow0}
		\frac{\log\Delta_{K(\cL)}(\delta)}
		{\log\delta} =1-\frac{\overline D}{3} ,\\
		\limsup_{n\to\infty}
		\frac{\log\Delta_{K(\cL)}(\delta_n)}
		{\log\delta_n}
		&=
		\limsup_{\delta\downarrow0}
		\frac{\log\Delta_{K(\cL)}(\delta)}{\log\delta} = 1-\frac{\underline D}{3}.
	\end{align*}
	Furthermore, $\frac{\log\delta_n}{\log 1/n}=\frac{\log v-\log n}{-\log n}\to 1$ as $n\to\infty$, which implies
	\begin{align*}
		\liminf_{n\to\infty}
		\frac{\log\Delta_{K(\cL)}(\delta_n)}{\log 1/n}
		&=
		1-\frac{\overline D}{3}, \qquad
		\limsup_{n\to\infty}
		\frac{\log\Delta_{K(\cL)}(\delta_n)}{\log 1/n}
		=
		1-\frac{\underline D}{3}.
	\end{align*}
	Finally, applying \eqref{eq:random-wet-part}
 	  proves the assertions concerning the expected missed area.
	
	For the expected number of vertices, Efron's identity from
	\cite{Efron1965}, see also \cite[Equation (8.12)]{SchneiderWeil}, says that
	\begin{equation}\label{eq:efron-identity}
		\mathbb E f_0(P_n)
		=
		\frac{n}{v}\,
		M_{n-1}.
	\end{equation}
	Replacing \(n\) by \(n-1\) does not change the lower or upper
	logarithmic order, since $\log(n-1)/\log n\to 1$ as $n\to\infty$.
	Taking logarithms in \eqref{eq:efron-identity} therefore yields
	\begin{align*}
		\liminf_{n\to\infty}
		\frac{\log\mathbb E f_0(P_n)}{\log n}
		&=1-\limsup_{n\to\infty}
		\frac{\log M_n}{\log 1/n}
		=1-\left(1-\frac{\underline D}{3}\right)
		=
		\frac{\underline D}{3},\\
		\limsup_{n\to\infty}
		\frac{\log\mathbb E f_0(P_n)}{\log n}
		&=1-\liminf_{n\to\infty}
		\frac{\log M_n}{\log 1/n}
        =		1-\left(1-\frac{\overline D}{3}\right)
		=
		\frac{\overline D}{3}.
	\end{align*}
	If \(D=\dim_M A_{\cL}\) exists, then
	\(\underline D=\overline D=D\), and the preceding limits are then
	equivalent to
	\[
	M_n
	=
	n^{-1+D/3+o(1)},\qquad \mathbb E f_0(P_n)
	=
	n^{D/3+o(1)}.
	\]
	The proof is thus complete.
\end{proof}

Under the assumptions of Corollary~\ref{cor:regular}, the estimates retain the slowly varying factors. Namely, if \(0<D\leq1\), then
\begin{equation*}
	\mathbb E\vol\bigl(K(\cL)\setminus P_n\bigr)
	\asymp
	n^{-1+D/3}L(n^{1/3})\qquad\text{and}\qquad \mathbb E f_0(P_n)
	\asymp
	n^{D/3}L(n^{1/3}).
\end{equation*}
In the critical case \(D=0\),
\[
\mathbb E\vol\bigl(K(\cL)\setminus P_n\bigr)
\asymp
\frac{\widehat L(n^{1/3})}{n}\qquad\text{and}\qquad
\mathbb E f_0(P_n)
\asymp
\widehat L(n^{1/3}).
\]
Here fixed positive factors in the arguments of \(L\) and
\(\widehat L\) have been suppressed. When $D=0$,
\(\widehat L\) is slowly varying by
\cite[Proposition~1.5.9a]{BinghamGoldieTeugels}.

\medspace

The preceding results show that the geometry of \(A_{\cL}\) determines
the asymptotic complexity of the associated random polygon in
\(K(\cL)\). It is therefore natural to compare these rates with the
two classical polygonal and smooth cases discussed in the introduction,
which go back to the foundational work of R\'enyi and Sulanke
\cite{RenyiSulanke1963,RenyiSulanke1964}.
If $K$ is a planar polygon, the expected missed area and the expected number of vertices satisfy
	\[
	\mathbb E\vol(K\setminus P_n)
	\asymp
	\frac{\log n}{n}
	\qquad\text{and}\qquad
	\mathbb E f_0(P_n)
	\asymp
	\log n,
	\]
	whereas, if \(K\) has a \(C^2\)-boundary with everywhere positive
	curvature, their respective orders are
	\[
	\mathbb E\vol(K\setminus P_n)\asymp n^{-2/3}
	\qquad\text{and}\qquad \mathbb E f_0(P_n)
	\asymp
	n^{1/3}.
	\]
For almost polygonal bodies associated with infinite fractal strings,
Remark~\ref{rem:strict-interpolation}, together with
\eqref{eq:random-wet-part} and Efron's identity
\eqref{eq:efron-identity}, shows that the expected missed area
and the expected number of vertices lie strictly between
their respective polygonal and smooth orders.

    The preceding estimates also determine the order of the fluctuations.
Indeed, the general variance estimates for planar random polygons
\cite[Corollary~1.4]{Pardon2012}, see also \cite{BaranySteiger2013},
together with \eqref{eq:random-wet-part}, imply
\begin{equation}\label{eq:VarEstimate}
    \operatorname{Var}(\vol(P_n))
    \asymp
    \frac{v}{n}\Delta_{K(\cL)}(v/n)
    \qquad\text{and}\qquad
    \operatorname{Var}(f_0(P_n))
    \asymp
    \frac{n}{v}\Delta_{K(\cL)}(v/n),
\end{equation}
where \(v=\vol(K(\cL))\). In particular, if
\(D=\dim_M A_{\cL}\) exists, these variances are of order
\(n^{-2+D/3+o(1)}\) and \(n^{D/3+o(1)}\), respectively.
The slowly varying refinements follow from
Corollary~\ref{cor:regular} in the same way.
Thus the geometry of the underlying fractal string also determines
the orders of the variances and hence the fluctuation scales in the
central limit theorems for planar random polygons established in
\cite[Corollary~1.2]{Pardon2012}.
	
For every \(\theta\in(0,1/3)\) and \(\gamma\in\R\),
Example~\ref{ex:logpower} provides a
fractal string \(\cL\)  with \(D=3\theta\) such that
\[
    \mathbb E\vol(K(\cL)\setminus P_n)
    \asymp n^{-1+\theta}(\log n)^{-\gamma}
    \qquad\text{and}\qquad
    \mathbb E f_0(P_n)
    \asymp n^\theta(\log n)^{-\gamma}.
\]
In particular, taking \(\gamma=0\) shows that every power strictly
between the classical polygonal and smooth powers can be realized.

The polygonal and smooth rates can be approached, but are not
attained by an infinite string. Towards the polygonal case, the
stretched-exponential strings of Example~\ref{ex:stretch} give,
for \(p>0\),
\[
    \mathbb E\vol(K(\cL)\setminus P_n)
    \asymp \frac{(\log n)^{1+1/p}}{n}
    \qquad\text{and}\qquad
    \mathbb E f_0(P_n)
    \asymp (\log n)^{1+1/p}.
\]
Towards the smooth case, Example~\ref{ex:iterated-log} shows that,
for every \(m\geq1\) and \(\gamma>0\), one can obtain
\[
    \mathbb E\vol(K(\cL)\setminus P_n)
    \asymp \frac{n^{-2/3}}{(\log_m n)^\gamma}
    \qquad\text{and}\qquad
    \mathbb E f_0(P_n)
    \asymp \frac{n^{1/3}}{(\log_m n)^\gamma}.
\]
Let us finally mention that the variance estimates \eqref{eq:VarEstimate}
show that all these asymptotic
behaviors also occur at the level of fluctuations. Indeed,
the variance of the number of vertices is comparable to its
expectation, while the variance of the area is comparable to
\(n^{-1}\) times the expected missed area. Thus the preceding
examples also determine the orders of the fluctuations in the
corresponding central limit theorems.

\subsection*{AI disclosure}

The connection between the asymptotic behavior of the convex floating body of almost polygonal convex bodies and the Minkowski dimension of the underlying sequence was first investigated in the Master's thesis of Hubertus Beyer \cite{Beyer2023}, supervised by the first author at Ruhr University Bochum. The second author subsequently suggested developing these ideas systematically in the language of fractal strings. The present paper builds on this earlier mathematical work.

During the preparation of the manuscript, the authors used ChatGPT (OpenAI) to assist with the detailed development and verification of mathematical arguments and with the drafting and a first version of the text. The underlying mathematical ideas, the research direction and also the applications to random polygons predate this use of AI. The authors have independently reviewed the mathematical content and take full responsibility for the results and their presentation.

\end{document}